\documentclass[11pt]{article}
 
\usepackage{graphicx}
\usepackage{color}
\usepackage{amsfonts}
\usepackage{amsmath}
\usepackage[mathcal]{eucal}
\usepackage{amsthm}
\usepackage{amssymb}
\usepackage{caption}
\usepackage{enumerate}
\usepackage{a4wide}
\usepackage{epstopdf}
\usepackage{mathtools}
\usepackage{soul}
\usepackage[normalem]{ulem}
\usepackage{float}
\usepackage{caption}

\newtheorem{theorem}{Theorem}[section]

\newtheorem{proposition}[theorem]{Proposition}
\newtheorem{corollary}[theorem]{Corollary}
\newtheorem{lemma}[theorem]{Lemma}

\newtheorem{conjecture}[theorem]{Conjecture}

\title{Clique immersions in graphs of independence number two with certain forbidden subgraphs}

\author{ Daniel A. Quiroz \\ \normalsize{Instituto de Ingenier\'ia Matem\'atica, Universidad de Valparaiso, Chile}  }

\date{}
\begin{document} 
\maketitle
\thispagestyle{empty}

\begin{abstract} 
The Lescure-Meyniel conjecture is the analogue of Hadwiger's conjecture for the immersion order. It states that every graph $G$ contains the complete graph $K_{\chi(G)}$ as an immersion, and like its minor-order counterpart it is open even for graphs with independence number~2. We show that every graph $G$ with independence number $\alpha(G)\ge 2$ and no hole of length between~$4$ and~$2\alpha(G)$ satisfies this conjecture. In particular, every $C_4$-free graph $G$ with $\alpha(G)= 2$ satisfies the Lescure-Meyniel conjecture. We give another generalisation of this corollary, as follows. Let $G$ and $H$ be graphs with independence number at most~2, such that $|V(H)|\le 4$. If $G$ is $H$-free, then $G$ satisfies the Lescure-Meyniel conjecture.



\end{abstract}

\section{Introduction}

For graphs $G$ and $H$, we say $G$ is \emph{$H$-free} if $G$ contains no induced subgraph isomorphic to~$H$. For a set of graphs $\mathcal{F}$, we say $G$ is \emph{$\mathcal{F}$-free} if $G$ is $F$-free for every $F\in \mathcal{F}$. A \emph{hole} in $G$ is an induced cycle of length at least four.  We let $\alpha(G)$ denote the \emph{independence number} (also known as \emph{stability number}) of~$G$. 

Hadwiger's conjecture~\cite{H43} asserts that every loopless graph $G$ contains the complete graph $K_{\chi(G)}$ as a minor. The conjecture is known to be true for $\chi(G)\le 6$, and  probably hard for larger values, as the proofs for cases $\chi(G)=5, 6$ already depend on the Four Colour Theorem. Nevertheless, Hadwiger's conjecture is known to be true for various graph classes. For instance, Reed and Seymour proved the conjecture for line graphs of multigraphs~\cite{RS04}. This was generalised  by Chudnovsky and Fradkin who proved it for quasi-line graphs~\cite{CF08}. 


One graph class for which Hadwiger's conjecture is still open is that of graphs with independence number at most 2. The importance of this case was perhaps first pointed out by Mader (see~\cite{PST03}), and is highlighted by Seymour in a recent survey~\cite{S16}. Plummer, Stiebitz and Toft~\cite{PST03} proved that for this class one can equivalently restate Hadwiger's conjecture as follows: every $n$-vertex graph $G$ with $\alpha(G)\le 2$ contains $K_{\lceil \frac n2 \rceil}$ as a minor. Previously, Duchet and Meyniel~\cite{DM82} had proved that every $n$-vertex graph $G$ contains a clique minor on $\lceil \frac{n}{2\alpha(G)-1}\rceil$ vertices. Despite much effort (see e.g.~\cite{B07, CS12, F10, FGS}), not much improvement has been obtained on this result, it still being open whether there is a constant $c>\frac 13$ such that every $n$-vertex graph with $\alpha(G)\le 2$ contains a clique minor on $\lceil cn\rceil$ vertices. However, some partial results are known in terms of forbidden subgraphs\footnote{For improvements on the theorem of Duchet and Meyniel for graphs with forbidden subgraphs and large independence number, see~\cite{BFS} and the references therein.}. Plummer et al.~\cite{PST03} proved that for every~$H$ on at most four vertices with $\alpha(H)\le 2$, if $G$ also has $\alpha(G)\le 2$ and is $H$-free then~$G$ satisfies Hadwiger's conjecture. Kriesell~\cite{K10} showed that this result still holds when~$H$ is allowed to have five vertices, while Bosse~\cite{B19}  showed it holds when $H$ is the wheel of six vertices. 


In this paper we focus on a conjecture related to Hadwiger's but concerning not minors but graph immersions. A graph~$G$ is said to contain another graph $H$ as an \emph{immersion} if there exists an injective function $\phi\colon V(H)\rightarrow V(G)$ such that:
\begin{enumerate}[(I)]
\item For every $uv\in E(H)$, there is a path in $G$, denoted $P_{uv}$, with endpoints $\phi(u)$ and~$\phi(v)$.
\item The paths in $\{P_{uv} \mid uv\in E(H) \}$ are pairwise edge disjoint.
\item The vertices of $\phi(V(H))$, called the \emph{branch vertices}, do not appear as interior vertices on paths~$P_{uv}$.
\end{enumerate}
Notice that if we further required the paths $P_{uv}$ to be internally vertex disjoint, instead of just edge disjoint as  in (II), then we would could say that $G$ contains $H$ as a topological minor. Thus if $G$ contains $H$ as a topological minor, then it contains $H$ as an immersion. On the other hand, the minor order and the immersion order are not comparable.


The immersion order has received a considerable amount of attention in recent years, particularly after Robertson and Seymour proved that it is a well-quasi-order~\cite{RS10}. Most of this attention has been directed towards the following conjecture of Lescure and Meyniel, which is an immersion-analogue of Hadwiger's conjecture.  (Abu-Khzam and Langston~\cite{AKL03} proposed a similar conjecture for a weaker notion of immersions requiring only (I) and (II).)

\begin{conjecture}[Lescure and Meyniel~\cite{LM89}]\label{hadimm}
Every graph $G$ contains $K_{\chi(G)}$ as an immersion.
\end{conjecture}

Whenever $\chi(G) \le 4$, this conjecture holds given the fact that Haj\'os' Conjecture 
 is true for these cases~\cite{D52}. The cases $5\le \chi(G) \le 7$  were established by Lescure and Meyniel~\cite{LM89} and by DeVos, Kawarabayashi, Mohar, and Okamura~\cite{Detal10}. But beyond this, there are few non-trivial cases for which this conjecture is known to hold. For instance, it is known that Conjecture~\ref{hadimm} holds for line graphs of simple graphs, but the conjecture is still open for line graphs of general multigraphs~\cite{GMcD19}.

For graphs with independence number at most 2, Vergara proposed the following conjecture which she showed is equivalent to Conjecture~\ref{hadimm} for this class of graphs. (A more general conjecture for graphs with arbitrary independence number is given in~\cite{BQSZ}.)
\begin{conjecture}[Vergara~\cite{V17}]\label{Vergara}
Every graph $G$ on $n$ vertices with $\alpha(G)\le 2$ contains an immersion of $K_{\lceil \frac n2\rceil}$.
\end{conjecture}
In support of her conjecture, Vergara showed that every $n$-vertex graph $G$ with $\alpha(G)\le 2$  contains $K_{\lceil \frac n3\rceil}$ as an immersion. Gauthier, Le and Wollan~\cite{GLW17} improved this as follows.
\begin{theorem}[Gauthier, Le and Wollan~\cite{GLW17}]\label{GLWdos}
Let $G$ be a graph on $n$ vertices with $\alpha(G)\le 2$. Then $G$ contains $K_{2\lfloor \frac n5 \rfloor}$ as an immersion.
\end{theorem}

One of our main results states that Conjecture~\ref{Vergara} holds for graphs with independence number~2 which exclude some subgraph on 4 or less vertices.  


\begin{theorem}\label{main2}
Let $G$ and $H$ be graphs with independence number at most 2, and $|V(H)|\le 4$. If $G$ is $H$-free, then $G$ satisfies Conjecture~\ref{Vergara}.
\end{theorem}

\begin{figure}[h]
 \centering
 \captionsetup{justification=centering}
 \bigskip
 \includegraphics[height=1.6 in]{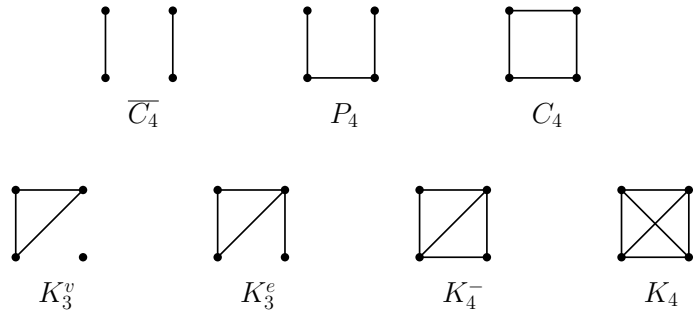} 
 \medskip
 \caption{The seven 4-vertex graphs with independence number at most 2.}
  \label{seven}
\end{figure}

Apart from providing evidence for Conjecture~\ref{Vergara}, Theorem~\ref{main2} is particularly interesting in light of the proof of Theorem~\ref{GLWdos}. That proof works by contradiction, assuming there is a counterexample~$G$   which is taken so as to minimise $|V(G)|+|E(G)|$. It is then shown that $G$ must contain an induced subgraph $H$ isomorphic to $C_5$. By choice of $G$, we have that $G-H$ satisfies Theorem~\ref{GLWdos}, and thus contains an immersion of a clique on $2\lfloor \frac {n-5}5 \rfloor$ vertices.  The proof finishes when it is shown that two vertices of $H$ can be added to the set of branch vertices of this immersion. Theorem~\ref{main2} implies that any minimal counterexample~$G$ to Conjecture~\ref{Vergara} must contain an induced copy of every 4-vertex graph $H$ with $\alpha(H)\le 2$ (Figure~\ref{seven} illustrates these graphs). If we could prove, that for some such $H$  we can add two vertices of $H$ to the immersion of a clique contained in $G-H$ (as was done for the~$C_5$), then we would prove Conjecture~\ref{Vergara}, i.e., Conjecture~\ref{hadimm} for graphs with independence number 2. We have, unfortunately, not been able to prove such a result.

Central to the proof of Theorem~\ref{main2} is showing that any $C_4$-free graph with $\alpha(G) = 2$ satisfies Conjecture~\ref{Vergara}. But we in fact prove a more general result as follows. An analogue statement was proved by Song and Thomas~\cite{ST17} for the so-called odd Hadwiger's conjecture. 

\begin{theorem}\label{forbholes}
Let $G$ be a graph with $\alpha(G)\ge 2$. If $G$ is $\{C_4,  \dots ,C_{2\alpha(G)}\}$-free, then $G$ satisfies Conjecture~\ref{hadimm}.
\end{theorem}

Improving this result to any graph $G$ with  no hole of length between $4$ and $2\alpha(G)-1$ (and $\alpha(G)\ge 3$) seems hard at the moment. Song and Thomas~\cite{ST17} proved that such graphs satisfy Hadwiger's conjecture, but for this they used that quasi-line graphs satisfy Hadwiger's conjecture. As mentioned earlier, Conjecture~\ref{hadimm} is open for line graphs of multigraphs and thus for quasi-line graphs.

As a tool for proving Theorem~\ref{forbholes}, we show that for every cycle  $C$ any \emph{inflation} (see Section~\ref{secholes} for the definition) of $C$ satisfies Conjecture~\ref{hadimm} (see Lemma~\ref{cycleinflation}). This is interesting in its own right, since Catlin~\cite{C79} showed that there are infinitely many inflations of odd cycles which are counterexamples to Haj\'os' conjecture (this includes graphs with independence number  2).

The rest of the paper is organised as follows. In Section~\ref{secholes} we prove Theorem~\ref{forbholes}. In Section~\ref{sechouse}, we prove that if $G$ with $\alpha(G)\le 2$ is $P_4$-free or $K_3^e$-free, then $G$ satisfies Conjecture~\ref{Vergara}. In Section~\ref{secowh} we do the same for graphs with $\alpha(G)\le 2$ which are $\overline{C_4}$-free or $K_3^v$-free. We finish the proof of Theorem~\ref{main2} in Section~\ref{secmain2}.

\section{Proof of Theorem~\ref{forbholes}}~\label{secholes}
In order to prove Theorem~\ref{forbholes} we first need a couple of lemmas, and, for these, some terminology and notation.  Set $\mathbb{N}_+:=\mathbb{N}\setminus\{0\}$. Let $G$ be a graph with $V(G)=\{v_1,v_2,\dots ,v_k\}$ and $f\colon V(G)\rightarrow\mathbb{N_+}$ a function. An \emph{$f$-inflation} of $G$ is a graph that can be obtained from $G$ by replacing each vertex $v_i\in V(G)$ with a clique with vertex set $B_i$ such that $|B_i|=f(v_i)$ and $B_i\cap B_j=\varnothing$ if $i\ne j$,  and by adding the edge $xy$, for $x\in B_i, y\in B_j$, if and only if $v_iv_j\in E(G)$.  If $H$ is an $f$-inflation of $G$ for some~$f$, then we say $H$ is an \emph{inflation} of $G$. For   $A, B \subseteq V(G)$ we say $A$ \emph{dominates} $B$ if every vertex in $B$ is adjacent to some vertex in~$A$. For ease of reading, we let $G\succcurlyeq H$ denote that $G$ contains $H$ as an immersion.

\begin{figure}[h]
 \centering
 \captionsetup{justification=centering}
 \bigskip
 \includegraphics[height=1.6 in]{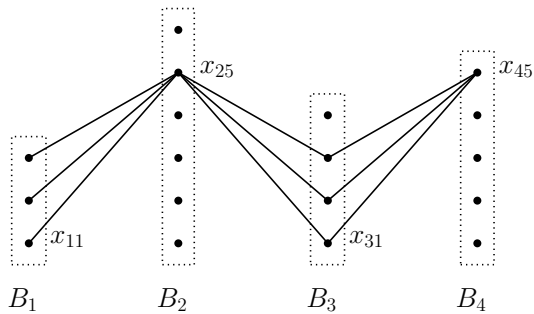} 
 \medskip
 \caption{The vertices of an inflation of a $P_4$. Only edges inside the paths $P_{r,5}$, $1\le r\le 3$, are shown. Only vertices in the path $P_{1,5}$ are labelled.}
  \label{pathinflation}
\end{figure}

\begin{lemma}\label{largepath}
Let $k\ge 1$ be an integer and $P=v_1,v_2, \dots, v_{2k}$ be a path. Let $f\colon V(P)\rightarrow \mathbb{N}_+$ be a function such that $p:=f(v_1)\le f(v_i)$ for every $i\in \{1,\dots ,2k\}$, and $q:=f(v_{2k})\le f(v_j)$ for every even $j\in \{2,\dots ,2k\}$. Let $G$ be an $f$-inflation of $P$. Then $B_1\cup B_{2k}$ are the branch vertices of an immersion of a clique on p+q vertices.
\end{lemma}
\begin{proof}
By definition of $f$, when $i$ is odd there are (distinct) vertices $x_{i1},\dots, x_{ip}$ in $B_i$, and when it is even there are vertices  $x_{i1},\dots, x_{iq}$ in $B_i$. Since $B_1$ and $B_{2k}$ both induce cliques, it suffices to show that there is a collection $\{P_{r,s}\mid 1\le r\le p,$ $1\le s\le q \}$ of mutually edge disjoint paths such that $P_{r,s}$ joins $x_{1r}$ with $x_{(2k)s}$, for every choice of $r$ and $s$. Let $P_{r,s}=x_{1r},x_{2s},x_{3r},x_{4s},\dots ,x_{(2k-1)r},x_{(2k)s}$. See Figure~\ref{pathinflation} for an illustration. Since $P_{r,s}$ is edge disjoint from $P_{r',s'}$ whenever $r\ne r'$ or $s\ne s'$, the result follows.
\end{proof}

As already mentioned, there are $f$-inflations of odd cycles (even for constant~$f$) which are counterexamples to Haj\'os' conjecture. We now show that Conjecture~\ref{hadimm} holds for inflations of cycles of any length.

\begin{lemma}\label{cycleinflation}
Let $k\ge 3$ be an integer. For a function $f\colon V(C_{k})\rightarrow \mathbb{N}_+$, let $G$ be an $f$-inflation of $C_{k}$. Then $G$ contains $K_{\chi(G)}$ as an immersion.
\end{lemma}
\begin{proof}
We proceed by induction on $k$, the result being trivially true when $k=3,4$. Let $k\ge 5$ and $v_1,\dots, v_k$  be the vertices of $C_k$ where $v_iv_{i+1}\in E(C_k)$ for $1\le i\le k-1$. Without loss of generality, assume that $|B_{k-1}|+|B_k|=\omega (G)$. We then have 
\begin{equation}\label{bagsize}
|B_1|\le |B_{k-1}| \mbox{ and } |B_{k-2}|\le |B_k|.
\end{equation}
Set $p:=\min \{|B_1|, |B_{k-2}|\}$, $q:=\max \{|B_1|, |B_{k-2}|\}$ and $P:=v_{k-2},v_{k-1},v_k,v_1$, and let $f'$ be the restriction of $f$ to $\{v_1,\dots, v_{k-2}\}$. Using Lemma~\ref{largepath} we see that $G$ contains an immersion~$G'$ which is an $f'$-inflation of the $C_{k-2}$ with vertices $v_1,\dots, v_{k-2}$, and with $v_iv_{i+1}\in E(C_{k-2})$ for $1\le i\le k-3$. By induction, we know that $G'\succcurlyeq K_{\chi(G')}$.

Let $c\colon (B_1\cup \dots \cup B_{k-2})\rightarrow \{1,\dots, \chi(G')\}$ be a proper colouring of $G'$ on $\chi(G')$ colours. We extend this colouring to a proper colouring of $G$ as follows. By \eqref{bagsize}, we can colour $|B_1|$ vertices of $B_{k-1}$ with $c(B_1)$, and $|B_{k-2}|$ vertices of $B_k$ with $c(B_{k-2})$. Notice that since $c$ is a proper colouring of $G'$ we have $c(B_1)\cap c(B_{k-2})=\varnothing$, which implies we have not created any monochromatic edge. If there remain uncoloured vertices in $B_{k-1}\cup B_k$ we use the colours in $\{1,\dots, \chi(G')\}\setminus c(B_1)\cup c(B_{k-2})$ to colour as many of them as possible. If there are no more uncoloured vertices in $B_{k-1}\cup B_k$, then we have properly coloured $G$ with $\chi(G')$ colours, and so $\chi(G)\le \chi(G')$. Hence, in this case we obtain $G\succcurlyeq G'\succcurlyeq K_{\chi(G')}\succcurlyeq K_{\chi(G)}$, as desired. If instead there remain uncoloured vertices, we colour them with the colours in $\{\chi(G')+1,\dots ,|B_{k-1}|+|B_k|\}$. This gives us a proper colouring of~$G$ with $|B_{k-1}|+|B_k|=\omega(G)$ colours. Therefore, $\chi(G)=\omega(G)$, and $G\succcurlyeq K_{\chi(G)}$ trivially. 
\end{proof}

We now prove Theorem~\ref{forbholes}. The proof uses Lemma~\ref{cycleinflation} after proving that if $G$ is \linebreak $\{C_4, \dots ,C_{2\alpha(G)}\}$-free and contains $C_{2\alpha(G)+1}$, then $V(G)$ can be partitioned into two sets~$A$ and $B$, such that $A$ induces an inflation of a cycle, while every vertex in $B$ is adjacent to all other vertices in~$V(G)$. For this, we follow~\cite[Theorem 2.2]{ST17}.

\begin{proof}[Proof of Theorem~\ref{forbholes}]
Notice that no graph $G$ can contain a hole of length $2\alpha(G)+2$ or larger. Hence, if $G$ is also $C_{2\alpha(G)+1}$-free then $G$ contains no hole, i.e., it is chordal. Since chordal graphs are perfect, we in particular have $\chi(G)=\omega(G)$, and $G\succcurlyeq K_{\chi(G)}$ trivially. 


Now we assume $G$ contains an induced copy $H$ of $C_{2\alpha(G)+1}$, on vertices $v_1,\dots ,v_{2\alpha(G)+1}$. We claim that every vertex in the graph $G - H$ (i.e. $G-V(H)$) is adjacent to either all of $V(H)$ or to exactly three consecutive vertices of~$H$. Let $u\in V(G- H)$ be such that it is not adjacent to some vertex of~$H$. Since $H$ contains an independent set of size $\alpha(G)$, $u$ must be adjacent to some vertex in~$H$ , so we assume $uv_1\notin E(G)$ while $uv_2\in E(G)$. Using that $G$ is $\{C_4, \dots ,C_{2\alpha(G)}\}$-free, it is not hard to see that $u$ cannot be adjacent to any of $v_5,\dots ,v_{2\alpha(G)+1}$. (Note this proves the claim, when $\alpha(G)=2$.) Since $\{u, v_6, v_8, \dots ,v_{2\alpha(G)}\}$ is an independent set of size $\alpha(G)$, then it  dominates the rest of the graph. Thus $uv_3, uv_4\in E(G)$, which gives the claim.

 For every $i\in \{1,\dots ,2\alpha (G)+1\}$, let $A_i$ be the set of vertices adjacent to $v_i,v_{i+1},v_{i+2}$ (all index operations in this proof are modulo $2\alpha (G)+1$), and  $A:=\bigcup_{i\in \{1,\dots ,2\alpha (G)+1\}}A_i\cup V(H)$. Notice that $G[A_i]$ is a clique, for if there were $u,w\in A_i$ with $uw\notin E(G)$, then $u,v_i,w,v_{i+2}$ would induce a~$C_4$. Also, for every $u\in A_i$ and $w\in A_{i+1}$ we must have $uw\in E(G)$, as otherwise $\{u,w,v_{i+4},v_{i+6},\dots ,v_{i+2\alpha (G)}\}$ would be an independent set of size $\alpha (G)+1$. Finally, notice that for every $u\in A_i$ and $w\in A_j$, with $j\notin\{i-1,i,i+1\}$, we have $uw\notin E(G)$ because $G$ is  $\{C_4, \dots ,C_{2\alpha(G)}\}$-free. Altogether, $A$ induces an inflation of $C_{2\alpha (G)+1}$.

Let $B$ be the set of vertices of $G$ which are adjacent to all of $V(H)$. Observe that $B$ induces a clique, for if there were $u,w\in B$ in with $uw\notin E(G)$ then $u,v_1,w,v_3$ would induce a~$C_4$.  Additionally, every vertex $u\in B$ is adjacent to every vertex $a\in A\setminus V(H)$, for otherwise $u,v_i,a,v_{i+2}$ would induce a $C_4$. Together, these two observations mean that $N(u)=V(G)-u$ for every $u\in B$. Therefore, we have $\chi(G[A])+|B|=\chi(G)$. Since $G[A]$ is an inflation of a cycle, Lemma~\ref{cycleinflation} tells us that $G[A]\succcurlyeq K_{\chi(G[A])}$, which implies that $G$ contains an immersion of a clique on  $\chi(G[A])+|B|=\chi(G)$ vertices, as desired.
\end{proof}

\begin{corollary}\label{C4}
Let $G$ be a $C_4$-free graph with $\alpha(G)\le 2$. Then $G$ satisfies Conjecture~\ref{Vergara}.
\end{corollary}

\section{Excluding the house}~\label{sechouse}
\begin{figure}[h]
 \centering
 \captionsetup{justification=centering}
 \bigskip
 \includegraphics[height=0.8 in]{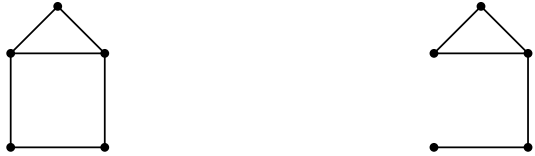}
 \medskip
 \caption{The \emph{house graph} and the \emph{one-wall-house graph}.}
  \label{houses}
\end{figure}

Let us say that a graph is \emph{house-free} if it excludes the \emph{house graph} given in Figure~\ref{houses}. In this section we show that every house-free graph $G$ with  $\alpha(G)\le 2$ satisfies Conjecture~\ref{Vergara}. The proof depends on  Corollary~\ref{C4} and the next lemma. Since the house graph contains both~$P_4$ and~$K_3^e$ as induced subgraphs, we obtain that every graph~$G$  with $\alpha(G)\le 2$ which  is $P_4$-free or $K_3^{e}$-free  satisfies Conjecture~\ref{Vergara}. 

A further piece of notation. For a graph $G$ and $v\in V(G)$, let $\bar N(v):=V(G)\setminus (N(v)\cup \{v\} )$.

\begin{lemma}\label{dominatingC4}
Let $G$ be an $n$-vertex graph with $\alpha(G)\le 2$. Suppose $G$ contains an induced copy~$H$ of~$C_4$, such that for every $e\in E(H)$ the  endpoints of $e$, together, dominate $G - H$, and $G-H$ satisfies Conjecture~\ref{Vergara}. Then $G$ contains $K_{\lceil \frac n2\rceil}$ as an immersion.
\end{lemma}
\begin{proof}
 Let $a_1,\dots, a_4$ be the vertices of~$H$, where $a_ia_{i+1}\in E(H)$ for $1\le i \le 3$. We remove $H$ from~$G$, and by hypothesis  $G-H$ contains an immersion of $K_{\lceil \frac {n-4}2\rceil}$ with branch vertices~$M$. Notice that it suffices to add two new branch vertices to this immersion. Only in one case will we be unable to do that, but then we will find some other clique (subgraph) on~$\lceil \frac n2\rceil$ vertices.

Set $Q:=V(G-H)\setminus M$. We clearly have,
\begin{equation}\label{QH5}
|Q|=n-4-\Big\lceil \frac {n-4}2\Big\rceil=\Big\lfloor \frac n2 \Big\rfloor -2.
\end{equation} 

Set $N_i:=N_G(a_i)$ and $\bar N_i:=\bar N_G(a_i)$ for every $i\in \{1,2,3,4\}$. We first show that we can assume 
\begin{equation}\label{nonneighMplus1C4}
|M\cap\bar N_i|\le |Q\cap N_i|+1, \mbox{ for every } i\in\{1,2,3,4\} 
\end{equation}
   
For suppose there was some $i\in \{1,2,3,4\}$ such that $|M\cap\bar N_i|-1> |Q\cap N_i|$. Then, by equality~\eqref{QH5}, we would have
\begin{equation*}
\Big\lfloor \frac n2 \Big\rfloor -2=|Q|< |Q\cap \bar N_i|+|M\cap\bar N_i|-1,
\end{equation*}
which implies $\big\lfloor \frac n2 \big\rfloor\le |\bar N_i- V(H)|$. Since $|\bar N_i\cap V(H)|= 1$ for any such $i$, we obtain $|\bar N_i|\ge \big\lfloor \frac n2 \big\rfloor +1\ge \big\lceil \frac n2 \big\rceil$. But given that $\alpha (G)\le 2$, we have that $\bar N_i$ induces a clique, which we just saw has at least~$\lceil \frac n2\rceil$ vertices. This proves~\eqref{nonneighMplus1C4}.

Since $\alpha(G)\le 2$, every set of non-adjacent vertices dominates the rest of the graph. Thus the condition that every edge $e$ of $H$ dominates $G - H$,  actually implies a stronger statement: every pair of vertices in $H$ dominates $G-H$. This means that 
\begin{equation}\label{disjointnonneigh}
\bar N_i \cap\bar N_j=\varnothing \mbox{ for every pair of distinct } i,j\in \{ 1,\dots ,4 \}.
\end{equation} 
By symmetry, we can assume without loss of generality that $|M\cap\bar N_1|\le |M\cap\bar N_i|$ for all $i\in \{2,3,4\}$. By~\eqref{disjointnonneigh} we then have, $|M\cap\bar N_1|\le \frac 13 (|M|-|M\cap\bar N_4|)$, which implies
 \begin{equation}\label{MellC4}
|M\cap\bar N_4|\le |M| - 3|M\cap\bar N_1|
\end{equation}

If $M\cap\bar N_4=\varnothing$, then our recent assumption also gives us $M\cap\bar N_1=\varnothing$. In this case, we trivially have an immersion of a clique with branch vertices $M\cup\{a_1,a_4\}$. So we can assume there is some vertex $z\in M\cap\bar N_4$. In light of~\eqref{nonneighMplus1C4}, we know there are injective functions mapping $M\cap\bar N_4\setminus\{z\}$ into $Q\cap N_4$. We choose one such injection $f$ such that $f(M\cap\bar N_4\setminus\{z\})\cap \bar N_1$ is as large as possible, and set $Q_f:=f(M\cap\bar N_4\setminus\{z\})$.

For every $x\in M\cap \bar N_4\setminus\{z\}$, consider the path $xa_i f(x)a_4$, for some $i  \in \{2,3\}$. This is possible since for every such $i$, $a_i$ dominates $M\cap \bar N_4$ and because $a_2,a_3$ together dominate $V(G- H)$. Since $f$ is an injection, these paths are mutually edge disjoint and so, together with the path $za_3a_4$, they attests that $G$ contains an immersion of a clique with branch vertices $M\cup \{a_4\}$. Our goal now is to show that the set of branch vertices can be expanded so as to also include $a_1$.

In fact, we shall now see that if we have 
\begin{equation}\label{listoC4}
|Q\cap N_1 \setminus Q_f|\ge |M\cap \bar N_1|-1,
\end{equation}
then we can add $a_1$ to the set of branch vertices in order to finish the proof. Choose a vertex $w\in M\cap \bar N_1$. If we have~\eqref{listoC4} then there is an injection $g\colon M\cap \bar N_1\setminus\{w\} \rightarrow Q\cap N_1 \setminus Q_f$.  For every $y\in M\cap \bar N_1\setminus\{w\}$ consider a path of the form $ya_i g(y)a_1$, for some $i\in \{2,3\}$. Using previous arguments it is easy  to see that these paths exists and are mutually edge disjoint. By definition, we have $\varnothing= g(M\cap \bar N_1\setminus\{w\})\cap Q_f= g(M\cap \bar N_1\setminus\{w\})\cap f(M\cap \bar N_4\setminus\{z\})$. This and~\eqref{disjointnonneigh} give us that the paths defined through $g$ are all edge disjoint from those defined through $f$. Then the paths $za_3a_4$, $wa_2a_1$, $a_1a_4$, and those defined through $f$ and $g$ are all mutually edge disjoint. Therefore all these paths, together, witness that $G$ contains an immersion of a clique with branch vertices $M\cup \{a_1, a_4\}$.

Now we need to check that \eqref{listoC4} indeed holds. Since $a_1$ and $a_4$ together dominate $V(G)-H$, we have $Q\cap \bar N_1\subseteq Q\cap N_4$. Then our choice of $f$ implies we either have $Q\cap\bar N_1\subseteq Q_f$ or $Q_f\subseteq Q\cap\bar N_1$. Therefore we obtain
\begin{equation}\label{lno1C4}
|Q\cap N_1 \setminus Q_f|= |Q\cap N_1|- \max\{ 0, |Q_f| - |Q\cap \bar N_1 |\}.
\end{equation}
When $\max\{ 0, |Q_f| - |Q\cap \bar N_1 |\}=0$, we use \eqref{nonneighMplus1C4} and \eqref{lno1C4} to see that we indeed have~\eqref{listoC4}. On the other hand, when $\max\{ 0, |Q_f| - |Q\cap \bar N_1 |\}\ne 0$ we start from \eqref{lno1C4} and actually obtain
\begin{align*}
|Q\cap N_1 \setminus Q_f|&\ge |Q|-  |Q_f| \\
&= |Q|-  |M\cap \bar N_4\setminus\{z\}| \\
&\ge |Q|- (|M| - 3|M\cap\bar N_1|)+1\\
&= 2\Big\lfloor \frac n2\Big\rfloor -n +3|M\cap\bar N_1| +1 \\
&\ge 3|M\cap\bar N_1|,
\end{align*}
where the second line comes from the definition of $f$, the third from \eqref{MellC4}, and the fourth from~\eqref{QH5}. The result follows.
\end{proof}

\begin{theorem}\label{house}
Let $G$ be a house-free graph with $\alpha(G)\le 2$. Then $G$ satisfies Conjecture~\ref{Vergara}. 
\end{theorem}
\begin{proof}
It suffices to show that every counterexample to Conjecture~\ref{Vergara} contains a copy of the house graph. For a contradiction, assume there is a house-free graph $G$ which is a counterexample to Conjecture~\ref{Vergara}. We choose $G$ with $|V(G)|$ minimum. By Corollary~\ref{C4},~$G$ contains an induced copy $F$ of $C_4$. Since $\alpha(G)\le 2$ every vertex in $G-F$ is adjacent to at least two consecutive vertices of $F$. But since $G$ is house-free we can further say that every such vertex must be adjacent to at least three consecutive vertices of $F$. Therefore, every edge of~$F$ satisfies that its  endpoints, together, dominate $G- F$. By minimality of $G$, $G-F$ satisfies Conjecture~\ref{Vergara}. However, by Lemma~\ref{dominatingC4}, this contradicts the choice of  $G$ as a counterexample to Conjecture~\ref{Vergara}. 
\end{proof}

\begin{corollary}~\label{P4}
Let $G$ be a graph with $\alpha(G)\le 2$. If $G$ is $P_4$-free or $K_3^e$-free, then $G$ satisfies Conjecture~\ref{Vergara}. 
\end{corollary}

\section{Excluding the one-wall house}\label{secowh}
We say a graph is \emph{o.w.h.-free} if it excludes the \emph{one-wall-house graph} given in Figure~\ref{houses}. In this section we prove that every o.w.h.-free graph $G$ with  $\alpha(G)\le 2$ satisfies Conjecture~\ref{Vergara}. The proof is slightly more complicated than that of Theorem~\ref{house}, depending on two lemmas similar to Lemma~\ref{dominatingC4}, and on Corollary~\ref{P4}. As a consequence we obtain that every graph~$G$  with $\alpha(G)\le 2$ which  is $\overline{C_4}$-free or $K_3^{v}$-free  satisfies Conjecture~\ref{Vergara}.  

\begin{lemma}\label{dominatingC5}
Let $G$ be an $n$-vertex graph with $\alpha(G)\le 2$. Suppose $G$ contains an induced copy~$H$ of $C_5$,  such that for every $e\in E(H)$ the  endpoints of $e$, together, dominate $G - H$, and every proper subgraph of $G$ satisfies Conjecture~\ref{Vergara}. Then $G$ contains $K_{\lceil \frac n2\rceil}$ as an immersion.
\end{lemma}

\begin{proof}
Let $a_1,\dots, a_5$ be the vertices of~$H$, where $a_ia_{i+1} \in E(H)$ for $1\le i \le 4$. We remove $P:=\{a_1,a_2,a_3,a_4\}$ from~$G$. By assumption we know that $G-P$ contains an immersion of $K_{\lceil \frac {n-4}2\rceil}$ with branch vertices~$M$. Thus, it suffices to add two new branch vertices to this immersion. 

Set $Q:=V(G)\setminus(P\cup M)$, and $N_i:=N_G(a_i)$, $\bar N_i:=\bar N_G(a_i)$ for every $i\in \{1,2,3,4\}$. As in the proof of Lemma~\ref{dominatingC4} we have
\begin{equation}\label{QH6}
|Q|=\Big\lfloor \frac n2 \Big\rfloor -2,
\end{equation} 
and can show we have
\begin{equation}\label{nonneighMplus1}
|M\cap\bar N_i|\le |Q\cap N_i|+1, \mbox{ for every } i\in\{1,2,3,4\}. 
\end{equation}
In fact, using that $|\bar N_i\cap P|=2$ when $i\in \{1,4\}$, we can further show that
\begin{equation}\label{nonneighM}
|M\cap\bar N_i|\le |Q\cap N_i|, \mbox{ for every } i\in\{1,4\}. 
\end{equation}

Also as in Lemma~\ref{dominatingC4}, the hypothesis implies the stronger statement that every pair of vertices in~$P$ dominates $G-H$. This means we have
\begin{equation}\label{disjointnonneighC5}
\bar N_i \cap\bar N_j \setminus\{ a_5\}=\varnothing \mbox{ for every pair of distinct } i,j\in \{ 1,\dots ,4 \}.
\end{equation}
 Let us fix $\ell$ such that $|M\cap\bar N_\ell\setminus\{ a_5\}|\le |M\cap\bar N_i\setminus\{ a_5\}|$ for all $i\in \{1,2,3\}$. By~\eqref{disjointnonneighC5} we have, $|M\cap\bar N_\ell\setminus\{ a_5\}|\le \frac 13 (|M\setminus\{ a_5\}|-|M\cap\bar N_4|)$, which implies
 \begin{equation}\label{Mell}
|M\cap\bar N_4|\le |M| - 3|M\cap\bar N_\ell\setminus\{ a_5\}|.
\end{equation}

By~\eqref{nonneighM}, there are injective functions mapping $M\cap\bar N_4$ into $Q\cap N_4$. We choose one such injection~$f$ such that $f(M\cap\bar N_4)\cap \bar N_\ell$ is as large as possible, and set $Q_f:=f(M\cap\bar N_4)$.

For every $x\in M\cap \bar N_4$, consider the path $xa_i f(x)a_4$, for some $i  \in \{1,2,3\}\setminus \{\ell\}$. This is possible since for every such $i$, $a_i$ dominates $M\cap \bar N_4$ and because for distinct $i,j  \in \{1,2,3\}\setminus \{\ell\}$, $a_i,a_j$ together dominate $G-H$. Since $f$ is an injection, these paths are mutually edge disjoint. This attests that $G$ contains an immersion of a clique with branch vertices $M\cup \{a_4\}$. Our goal now is to show that the set of branch vertices can be expanded so as to also include $a_\ell$.

Towards this new goal, we first note that we can assume
\begin{equation}\label{detalle}
|M\cap \bar N_\ell|\ge 1, \mbox{ and \,} |M\cap \bar N_\ell|\ge 2 \mbox{ when } \ell\ne 1.
\end{equation}
For otherwise, if $|M\cap \bar N_\ell|=0$, the path joining $a_\ell$ and $a_4$ in  $G[P]$, together with those defined through~$f$, guarantees that $G$ contains an immersion of a clique with branch vertices $M\cup \{a_\ell, a_4\}$, as desired. If instead we have  $\ell\ne 1$ and $M\cap \bar N_\ell=\{z\}$ for some vertex $z$,  we let $P_1, P_4$ be the paths joining $a_\ell$ with $a_1$ and $a_4$, respectively, in $G[P]$. Since $a_1,a_\ell$ together dominate~$M$, then $a_1$ is adjacent to $z$. Therefore, the paths $P_1z$, $P_4$ and those defined through $f$ witness that $G$ has an immersion of a clique with branch vertices $M\cup \{a_\ell, a_4\}$. This proves~\eqref{detalle}.

It is not hard to see that if we have
\begin{equation}\label{listo}
|Q\cap N_\ell \setminus (Q_f\cup\{ a_5\})|\ge |M\cap \bar N_\ell|,
\end{equation}
then we can add $a_\ell$ to the set of branch vertices in order to finish the proof.  From here onwards, the proof requires some case analysis. In some cases we will be able to prove \eqref{listo} holds, and thus obtain the desired clique immersion. In others, we will be one vertex away from having \eqref{listo} and will require some more care.

We have that $a_\ell$ and $a_4$ dominate $V(G)\setminus P$. Thus $Q\cap \bar N_\ell\subseteq Q\cap N_4$, which by our choice of $f$ implies we either have $Q\cap\bar N_\ell\subseteq Q_f$ or $Q_f\subseteq Q\cap\bar N_\ell$. Given that when $\ell\ne 1$ we have $a_\ell a_5\notin E(G)$, for these values of $\ell$ we obtain
\begin{equation}\label{lno1}
|Q\cap N_\ell \setminus (Q_f\cup\{ a_5\})|= |Q\cap N_\ell|- \max\{ 0, |Q_f| - |Q\cap \bar N_\ell |\}.
\end{equation}
This also holds when $\ell =1$ and $a_5\in M$, while when $\ell =1$ and $a_5\in Q$ we only have
\begin{equation}\label{lsi1}
|Q\cap N_\ell \setminus (Q_f\cup\{ a_5\})|\ge |Q\cap N_\ell|- \max\{ 0, |Q_f| - |Q\cap \bar N_\ell |\}-1. 
\end{equation}

We first deal with the case $\ell=1$, and consider three subcases. When $\ell=1$ and $\max\{ 0, |Q_f| - |Q\cap \bar N_\ell |\}\ne 0$ we start from \eqref{lsi1} in order to see that
\begin{align*}
|Q\cap N_\ell \setminus (Q_f\cup\{ a_5\})|&\ge |Q|-  |Q_f| -1\\
&= |Q|-  |M\cap \bar N_4| -1\\
&\ge |Q|- (|M| - 3|M\cap\bar N_\ell|)-1\\
&= 2\Big\lfloor \frac n2\Big\rfloor -n +3|M\cap\bar N_\ell| -1\\
&\ge 3|M\cap\bar N_\ell|-2,
\end{align*}
where the second line comes from the definition of $f$; the third from \eqref{Mell} and the fact that $|M\cap\bar N_\ell\setminus\{ a_5\}|=|M\cap\bar N_\ell|$ when $\ell=1$; and the fourth from \eqref{QH6}. In light of \eqref{detalle}, for this subcase we obtain \eqref{listo}, which guarantees the desired immersion.

When $\ell=1$, $\max\{ 0, |Q_f| - |Q\cap \bar N_\ell |\}= 0$, and $a_5\in M$ we easily obtain \eqref{listo} by combining~\eqref{nonneighM} and~\eqref{lno1}.

In the remaining subcase, we combine~\eqref{nonneighM} and~\eqref{lsi1} to see that $|Q\cap N_1 \setminus (Q_f\cup\{ a_5\})|\ge |M\cap \bar N_1|-1$. We then chose a vertex $z\in M\cap \bar N_1$, and an injection $g\colon M\cap \bar N_1\setminus\{z\}\rightarrow Q\cap N_1 \setminus (Q_f\cup\{ a_5\})$. For every $y\in M\cap \bar N_1\setminus\{z\}$ consider a path of the form $xa_ig(x)a_1$, for $a_i$ in $\{2,3\}$. As in the proof of Lemma~\ref{dominatingC4}, it is not hard to show that such paths exists, are mutually edge disjoint and disjoint from all the paths defined using $f$. To see that $G$ contains an immersion of a clique with branch vertices $M\cup \{a_1, a_4\}$, it remains to show that we can find  mutually edge disjoint paths joining $a_1$ to $z$ and to~$a_4$, and which are also edge disjoint from all previously considered paths. For this, we take the paths $a_1a_2z$ and $a_1a_5a_4$. The first path is edge disjoint from the paths defined using $g$ for none of them contain $z$. By~\eqref{disjointnonneighC5} the paths defined using $f$ do not contain~$z$ either, and are thus edge disjoint from $a_1a_2z$, as well. By definition, the paths defined using $g$ are also edge disjoint from $a_1a_5a_4$. Finally, as noted before, the definition of $f$ tells us that $\max\{ 0, |Q_f| - |Q\cap \bar N_1 |\}= 0$ implies $Q_f\subseteq Q\cap\bar N_1$. Therefore, we must have $a_5\notin Q_f$, since $a_5\in N_1$ (indeed, in this subcase, $a_5\in Q\cap N_1$). This witnesses that the paths defined using~$f$ are edge disjoint from the path $a_1a_5a_4$. Altogether, we obtain the desired  immersion with branch vertices $M\cup \{a_1, a_4\}$.

To conclude we now consider the cases when $\ell\ne 1$, and for that two subcases. In none of the subcases do we obtain~\eqref{listo}, and only obtain 
\begin{equation}\label{casilisto}
|Q\cap N_\ell \setminus (Q_f\cup\{ a_5\})|\ge |M\cap \bar N_\ell|-1.
\end{equation}
However, it is not hard to see that, for this value of $\ell$, this suffices.

When $\max\{ 0, |Q_f| - |Q\cap \bar N_\ell |\}\ne 0$, we proceed similarly as when we had $\ell=1$ and $\max\{ 0, |Q_f| - |Q\cap \bar N_\ell |\}\ne 0$, but start from \eqref{lno1} so as to obtain 
$$|Q\cap N_\ell \setminus (Q_f\cup\{ a_5\})|\ge 3|M\cap\bar N_\ell\setminus\{ a_5\}|-1,$$ which in light of \eqref{detalle} indeed gives  \eqref{casilisto}. When $\max\{ 0, |Q_f| - |Q\cap \bar N_\ell |\}=0$, we instead use \eqref{nonneighMplus1} and \eqref{lno1}, to see that again we have~\eqref{casilisto}. The result follows.
\end{proof}

\begin{lemma}\label{dominatingP4}
Let $G$ be an $n$-vertex graph with $\alpha(G)\le 2$. Suppose $G$ contains an induced copy~$H$ of $P_4$,  such that for every  $e\in E(H)$ the  endpoints of $e$, together, dominate $G - H$, and $G-H$ satisfies Conjecture~\ref{Vergara}. Then $G$ contains $K_{\lceil \frac n2\rceil}$ as an immersion.
\end{lemma}

The proof of Lemma~\ref{dominatingP4} is so similar to the proof of Lemma~\ref{dominatingC5} that we do not include it. In the proof of Lemma~\ref{dominatingC5} we take out a copy $P$ of $P_4$ from $G$ and add two of its vertices to the immersion guaranteed, by induction,  in $G-P$. To prove Lemma~\ref{dominatingP4} we can proceed in the same way. The proof is in fact much easier, since we do not need to deal with the vertex $a_5$ of the proof of Lemma~\ref{dominatingC5}. In particular, this means that in every case we obtain~\eqref{lno1}, and never only~\eqref{lsi1}. Hence, when $\ell=1$ we always obtain~\eqref{listo}, which allows us to avoid the most delicate case of the proof of Lemma~\ref{dominatingC5}.

\begin{theorem}\label{owh}
Let $G$ be a o.w.h.-free graph with $\alpha(G)\le 2$. Then $G$ satisfies Conjecture~\ref{Vergara}. 
\end{theorem}
\begin{proof}
It suffices to show that every counterexample to Conjecture~\ref{Vergara} contains a copy of the one-wall-house graph. By Corollary~\ref{P4}, every counterexample $G$ to Conjecture~\ref{Vergara} contains a copy $P$ of $P_4$. We choose $G$ with $|V(G)|$ minimum, that is, every proper subgraph of $G$ satisfies Conjecture~\ref{Vergara}. We let $P=a_1,a_2,a_3,a_4$. Since $\alpha(G)\le 2$, every vertex of $G-P$ is adjacent to at least two consecutive (modulo 4) vertices of $P$. But we must have some vertex in $G-P$ being adjacent to exactly two vertices of $P$. Otherwise, every edge of $P$ would dominate $G-P$, and by Lemma~\ref{dominatingP4} this would contradict our choice of $G$. Since every vertex of $G-P$ must also be adjacent to one of $a_1$ and~$a_4$, we conclude that $G$  contains a copy of the one-wall-house graph or a copy of $C_5$.

We assume $G$ contains a copy $H$ of $C_5$, with vertices $a_1,\dots ,a_5$. Since $\alpha(G)\le 2$ every vertex of $G-H$  is adjacent to at least three consecutive  vertices of $H$. But some vertex $v$ in $G-H$ must be adjacent to exactly three vertices, $a_1,a_2,a_3$ say, of $H$. Otherwise, every edge of $H$ would dominate $G-H$, which by Lemma~\ref{dominatingC5} would contradict the choice of $G$. Then, $v,a_2,a_3,a_4,a_5$ induce a copy of the one-wall-house graph.
\end{proof}

\begin{corollary}~\label{overlineC4}
Let $G$ be a graph with $\alpha(G)\le 2$. If $G$ is $\overline{C_4}$-free or $K_3^v$-free, then $G$ satisfies Conjecture~\ref{Vergara}. 
\end{corollary}

\section{Proof of Theorem~\ref{main2}}\label{secmain2}
Corollaries~\ref{C4}, \ref{P4} and \ref{overlineC4}, and the next two propositions immediately imply Theorem~\ref{main2}. Notice that while Corollary~\ref{overlineC4} depends on Corollary~\ref{P4}, and this one on Corollary~\ref{C4}, the following propositions are independent of all our previous results. Notice also that for $K_4^-$-free graphs we actually find a clique subgraph on $\lceil \frac n2 \rceil$ vertices.

\begin{proposition}\label{K4}
Let $G$ be a $K_4$-free $n$-vertex graph with $\alpha(G)\le 2$. Then $G$ satisfies Conjecture~\ref{Vergara}.
\end{proposition}
\begin{proof}
The result is trivially true when $n\le 4$, while by a result from Ramsey Theory~\cite{GG55}, if $G$ is a $K_4$-free graph with $\alpha(G)\le 2$, then $n\le 8$.  Consider first the case $n=5$. Since $\alpha(G)\le 2$, for every vertex $v$, $\bar N(v)$ induces a clique. Thus we  assume $|\bar N(v)|\le 2$, for every $v\in V(G)$. Let $v_1\in V(G)$ be adjacent to vertices $v_2,v_3$. If $v_2v_3\in E(G)$ or $N(v_2)\cap N(v_3)\ne \{v_1\}$, then $G$ contains a cycle, and trivially an immersion of~$K_3$. So we assume $v_2v_4,v_3v_5\in E(G)$. Since $|N(v_4)|\ge 2$, we inevitably have a cycle and thus an immersion of~$K_3$. The case $n=5$ follows.

The case $n=6$ follows trivially from the case $n=5$. But more can be said. A well-known result of Ramsey Theory gives that if $\alpha(G)\le 2$ and $n\ge 6$, then~$G$ contains a $K_3$ subgraph. Consider now the case $n=7$. By the result just mentioned $G$ contains a clique subgraph with set of vertices $M=\{m_1, m_2, m_3\}$. The graph $G-M$ on vertices $a_1,a_2,a_3,a_4$ must contain a set of non-adjacent vertices, say $a_1,a_2$, or else it would be isomorphic to a $K_4$. Also none of $a_1,a_2$ can be adjacent to all of $M$. Since $a_1,a_2$ are non-adjacent they dominate the rest of the graph, so we can assume without loss of generality that $a_1m_1, a_2m_2, a_2m_3\in E(G)$, while $a_1m_3,a_2m_1\notin E(G)$. Since $\alpha(G)\le 2$, $\bar N(v)$ is a clique for every $v\in V(G)$, and since~$G$ is $K_4$-free we have $|\bar N(v)|\le 3$ for every $v$. This, together with our assumption that $a_1m_3,a_2m_1\notin E(G)$ gives us that each of $a_1,a_2$ must be adjacent to one of $a_3,a_4$. If they are both  adjacent to the same vertex $a_i$, $i\in \{3,4\}$, then we can find an immersion of a clique on $\lceil \frac n2 \rceil$ vertices, $a_2,m_1,m_2,m_3$, by using the path $a_2,a_i,a_1,m_1$. So we can assume $a_1a_3,a_2a_4\in E(G)$, while $a_1a_4,a_2a_3\notin E(G)$. Since $|\bar N(a_1)|\le 3$ we must have $a_1m_2\in E(G)$. This completely specifies the neighbourhoods of $a_1$ and $a_2$. In particular, each has exactly two neighbours in $M$.

If we have $a_3a_4\in E(G)$, then we can obtain a clique immersion with branch vertices $M\cup\{a_2\}$ by using the path $a_2a_4a_3a_1m_1$. So we assume we do not have this edge. Then, by symmetry, the arguments given in the previous paragraph, tell us that $a_3$ and $a_4$ each have exactly two neighbours in $M$. Note that we cannot have $N(a_4)\cap M=\{m_2,m_3\}$, for this would imply that $m_2,m_3,a_2,a_4$ induce $K_4$ subgraph.   Thus we have $a_4m_1\in E(G)$, and we can obtain the desired immersion using the path $a_2a_4m_1$. This finishes the case $n=7$.

The case $n=8$ follows trivially from the case $n=7$, and the result follows.
\end{proof}

\begin{proposition}\label{K4-}
Let $G$ be a $K_4^-$-free $n$-vertex graph with $\alpha(G)\le 2$. Then $G$ contains $K_{\lceil \frac n2 \rceil}$ as a subgraph.
\end{proposition}
\begin{proof}
We show that every counterexample $G$ to Conjecture~\ref{Vergara} contains a copy of~$K_4^-$. If $G$ contains a copy $H$ of $C_5$, then $\alpha(G)\le 2$ implies that every vertex of $G-H$  is adjacent to at least three consecutive  vertices of $H$. Since this gives that $G$ contains a copy of $K_4^-$, we assume that~$G$ is $C_5$-free.

We may further assume that $G$ is not a clique and that $n\ge 5$. Let $x$ be a vertex of minimum degree. Then $\bar N(x)$ is a clique on at least one vertex, and we can assume $N(x)$ is not a clique. Since $G$ is $K_4^-$-free, $N(x)$ cannot contain an induced path on three vertices. It follows that $N(x)$ is the disjoint union of two cliques, say $A$ and $B$. Since we assumed $N(x)$ is not a clique, neither of these cliques can be empty.

We claim that either $A$ or $B$ is complete to $\bar N(x)$. Suppose there are $a\in A$, $b\in B$ and $y,z\in \bar N(x)$ such that $ay, bz\notin E(G)$. Then $by, az\in E(G)$ and so $y\ne z$. But then, $G[\{x,a,z,y,b\}]$ is a $C_5$, a contradiction. Thus, either $A$ or $B$, say $A$ is complete to $\bar N(x)$, as claimed. Then $V(G)$ is the disjoint union of two cliques $A\cup \bar N(x)$ and $B\cup\{x\}$. The result follows. 
\end{proof}

\section*{Acknowledgements}
The author thankfully acknowledges support from FONDECYT/ANID Iniciaci\'on en Investigaci\'on Grant 11201251; from \mbox{ANID} + PIA/Concurso Apoyo a Centros Cient\'ificos
y Tecnol\'ogicos de Excelencia con Financiamiento Basal, AFB170001; from Programa Regional MATH-AMSUD, MATH190013; and from Concurso para Proyectos de Investigaci\'on Conjunta ANID Chile -- FAPESP Brasil, 2019/13364-7.

The author wishes to thank an anonymous referee for finding and fixing some minor errors in a preliminary version, and for providing the present, shorter, proof of Proposition~\ref{K4-}.

\end{document}